\documentclass[12pt]{amsart} 
\usepackage[mathscr]{eucal}
\usepackage{amsmath,amsfonts}
\parskip=\smallskipamount

\newtheorem{theorem}{Theorem}[section]

\newtheorem{proposition}[theorem]{Proposition}
\newtheorem{corollary}[theorem]{Corollary}
\newtheorem{lemma}[theorem]{Lemma}

\newcommand{\CC}{{\mathbb C}}
\newcommand{\NN}{{\mathbb N}}

\newcommand{\ZZ}{{\mathbb Z}}
\newcommand{\DD}{{\mathbb D}}
\newcommand{\RR}{{\mathbb R}}
\newcommand{\FF}{{\mathbb F}}
\newcommand{\TT}{{\mathbb T}}

\makeatletter
\@addtoreset{equation}{section}
\makeatother

\newcommand{\cA}{{\mathcal A}}
\newcommand{\cB}{{\mathcal B}}
\newcommand{\cC}{{\mathcal C}}
\newcommand{\cD}{{\mathcal D}}
\newcommand{\cE}{{\mathcal E}}
\newcommand{\cF}{{\mathcal F}}
\newcommand{\cG}{{\mathcal G}}
\newcommand{\cH}{{\mathcal H}}
\newcommand{\cK}{{\mathcal K}}
\newcommand{\cL}{{\mathcal L}}
\newcommand{\cM}{{\mathcal M}}
\newcommand{\cN}{{\mathcal N}}

\newcommand{\cP}{{\mathcal P}}
\newcommand{\cR}{{\mathcal R}}

\newcommand{\cT}{{\mathcal T}}
\newcommand{\cJ}{{\mathcal J}}

\newdimen\expt
\def\boxit#1{\setbox0\hbox{$\displaystyle{#1}$}
      \hbox{\lower.4\expt
 \hbox{\lower3\expt\hbox{\lower\dp0
      \hbox{\vbox{\hrule height.4\expt
 \hbox{\vrule width.4\expt\hskip3\expt
      \vbox{\vskip3\expt\box0\vskip2\expt}%
 \hskip3\expt\vrule width.4\expt}\hrule height.4\expt}}}}}}

\newcommand{\be}{\begin{equation}}                         
\newcommand{\ee}{\end{equation}}
\newcommand{\bs}{\begin{subequations}}                         
\newcommand{\es}{\end{subequations}}
\newcommand{\ba}{\left[ \begin{array}}
\newcommand{\ea}{\\ \end{array} \right]}

\begin{document}

\title %[] 
{Tensor Algebras and Displacement Structure. 
III. Asymptotic properties} 

\author{M.~Barakat} 
\author{T.~Constantinescu}

\address{}

\address{Department of Mathematics,
  University of Texas at Dallas,
  Box 830688 Richardson, TX 75083-0688, U. S. A.}
\email{tiberiu@utdallas.edu}
\email{mxb6460@dcccd.edu}

\maketitle

\noindent
{\small {\bf Astract.}
We continue to investigate some classes of Szeg\"o type
polynomials in several variables. We focus on asymptotic properties of these 
polynomials and we extend several classical results of G. ~Szeg\"o to this 
setting.
}


\noindent
{\small {\bf AMS subject classicafication:} 15A69, 47A57}

\section{Introduction}
An extension to several non-commuting variables of the Szeg\"o 
orthogonal 
polynomials on the unit circle was considered in \cite{CJ}.
Some of the basic algebraic results on these polynomials were also obtained, 
including recurrence equations, Christoffel-Darboux formulae, and a Favard type
result. Also, it was explained their connection with displacement 
structure theory. Our main goal is to continue to investigate this kind of 
polynomials
and in this paper we focus on some of their asymptotic properties. 
There are several fundamental results of G. Szeg\"o involving asymptotic
properties of the orthogonal polynomials on the unit
circle. Thus, let $\TT $ be the unit circle and 
let $\mu $ be a positive Borel measure on $\TT $ with $\log \mu '\in L^1$.
Also let $\{\varphi _n\}_{n\geq 0}$ be the family of orthogonal 
polynomials associated to $\mu $ and $\varphi _n^{\sharp }(z)=z^n
\overline{\varphi _n(1/\overline{z})}$, $n\geq 0$. 
It is well-known (see \cite{Sz}) that 
\begin{equation}\label{unu}
\varphi _n\rightarrow 0
\end{equation}
and 
\begin{equation}\label{doi}
\frac{1}{\varphi _n^{\sharp }}\rightarrow \Theta _{\mu },
\end{equation}
where $\Theta _{\mu}$ is the spectral factor of $\mu $ and the convergence is 
uniform on the compact subsets of the unit disk $\DD $. The second limit 
\eqref{doi} is related to the so-called Szeg\"o limit theorems concerning
the asymptotic behaviour of Toeplitz determinants. Thus, 
$$\frac{\det T_n}{\det T_{n-1}}=\frac{1}{|\varphi _n^{\sharp }(0)|^2},
$$
where $T_n=[s_{i-j}]_{i,j=0}^n$ and $\{s_k\}_{k\in \ZZ }$ is the set
of the Fourier coefficients of $\mu $. As a consequence of the previous 
relation and \eqref{doi} we deduce Szeg\"o's first limit theorem,
\begin{equation}\label{trei}
\lim _{n\rightarrow \infty }
\frac{\det T_n}{\det T_{n-1}}=|\Theta _{\mu }(0)|^2=
\exp (\frac{1}{2\pi }\int _0^{2\pi }\log \mu '(t)dt).
\end{equation}
The second (strong) Szeg\"o limit theorem improves \eqref{trei}
by showing that 
\begin{equation}\label{patru}
\lim _{n\rightarrow \infty }
\frac{\det T_n}{g^{n+1}(\mu )}=
\exp \left(\frac{1}{\pi }\int\int _{|z|\leq 1}
|\Theta '_{\mu }(z)/\Theta _{\mu }(z)|^2d\sigma (z)\right),
\end{equation}
where $g(\mu )$ is the limit in \eqref{trei} and $\sigma $ is the 
planar Lebesgue measure. These two limits \eqref{trei}
and \eqref{patru} have an useful interpretation in terms of asymptotics 
of angles in the geometry of a stochastic process associated to $\mu $
(see \cite{GS}). 

Our goal in this paper is to extend these results to the class of orthogonal
polynomials on several non-commuting variables introduced
in \cite{CJ}. 
The paper is organized as follows. In Section~2 we review notation
and a framework for studying orthogonal polynomials associated to 
polynomial relations on several non-commuting variables. 
Thus the paper can be read independently of \cite{Co1} and 
\cite{CJ}. In Section ~3 
we analyse the case of no relation in dimension one.
It turns out that this is, in fact, the most general situation, and for this
reason we treat this case separately. The main result is Theorem ~3.3, 
which extends \eqref{unu} and \eqref{doi}. Theorem ~3.4 contains extensions
of \eqref{trei} and \eqref{patru}. In Section ~4 we discuss
a few examples. First, we show how to recapture the 
classical setting of orthogonal polynomials on the unit circle and on the real
line. Then, we turn our attention to the orthogonal polynomials 
considered in \cite{CJ}.

\section{Preliminaries}
We introduce some necessary
terminology and notation. Especially, we briefly review a rather 
familiar setting for orthogonal 
polynomials associated to relations on several variables
(for some details, see \cite{Co}). 

\subsection{Tensor Algebras}
Let $\FF _N^+$ be the unital free semigroup 
on $N$ generators $1,\ldots ,N$ with lexicographic
order $\prec $. The empty word is the identity element
and the length of the word $\sigma $ is 
denoted by $|\sigma |$. The length of the empty word is $0$
and $l(\sigma )$ denotes the number of words $\tau \preceq\sigma $.

The tensor algebra over $\CC ^N$ 
is defined by the algebraic direct sum 
$$\cT _N=\oplus _{k\geq 0}(\CC ^N)^{\otimes k},$$
where $(\CC ^N)^{\otimes k}$ denotes the $k$-fold tensor
product of $\CC ^N$ with itself. The addition is the 
componentwise addition and the  multiplication is defined
by juxtaposition:
$$(x\otimes y)_n=\sum _{k+l=n}x_k\otimes y_l.$$
If $\{e_1,\ldots ,e_N\}$ is the standard basis of
$\CC ^N$, then $\{e_{i_1}\otimes \ldots \otimes e_{i_k}
\mid 1\leq i_1, \ldots ,i_k\leq N\}$ is an orthonormal 
basis of $\cT _N$. 
If $\sigma =i_1\ldots i_k$ then we write
$e_{\sigma }$ instead of $e_{i_1}\otimes \ldots \otimes e_{i_k}$, 
so that any element of $\cT _N$ can be uniquely written in the 
form
$x=\sum _{\sigma \in \FF _N^+}c_{\sigma }e_{\sigma }$,
where only finitely many of the complex numbers $c_{\sigma }$ 
are different from $0$.

Another construction of $\cT _N$ is given by the 
algebra $\cP _N$ of polynomials in $N$ noncommuting
indeterminates $X_1,\ldots ,X_N$ with complex coefficients.
Each element $P\in \cP _N$ can be uniquely written 
in the form $P=\sum _{\sigma \in \FF _N^+}c_{\sigma }X_{\sigma }$,
with $c_{\sigma }\ne 0$ for finitely many $\sigma $'s and 
$X_{\sigma }=X_{i_1}\ldots X_{i_k}$ where $\sigma =i_1\ldots i_k
\in \FF _N^+$. The linear extension $\Phi _1$ of the mapping 
$e_{\sigma }\rightarrow X_{\sigma }$, $\sigma \in \FF _N^+$, 
gives an isomorphism of $\cT _N$ with $\cP _N$.

Another known realization of the tensor algebra was used in \cite{CJ}
in order to establish a connection with the displacement structure
theory. This was useful since many results for the tensor algebra 
could be seen
just as particular instances of more general results in the triangular
algebra. Thus, let $\cE $ be a Hilbert space
and define: $\cE _0=\cE $ and for $k\geq 1$, 
\begin{equation}\label{prebaza}
\cE _k=\underbrace{\cE _{k-1}\oplus \ldots \oplus \cE _{k-1}}_
{N \,terms}=\cE ^{\oplus N}_{k-1}.
\end{equation}
For $\cE=\CC $ we have that $\CC _k$
can be identified with $(\CC ^N)^{\otimes k}$ and $\cT _N$ is
isomorphic to the algebra $\cL _N$ of lower triangular 
operators $T=[T_{k,j}]\in \cL (\oplus _{k\geq 0}\CC _k)$
with the property 
\begin{equation}\label{abaza}
T_{k,j}=\underbrace{T_{k-1,j-1}\oplus \ldots \oplus T_{k-1,j-1}}
_{N \, terms}=T^{\oplus N}_{k-1,j-1},
\end{equation}
for $k\leq j$, $k,j\geq 1$, and $T_{j,0}=0$ for all 
sufficiently large $j's$. The isomorphism is given by the 
map $\Phi _2$  defined as follows:
let $x=(x_0,x_1,\ldots )\in \cT _N$ ($x_p\in (\CC ^N)^{\otimes p}$
is the $p$th homogeneous component of $x$); then 
$x_p=\sum _{|\sigma |=p}c_{\sigma }e_{\sigma }$
and for $j\geq 0$, $T_{j,0}$ is given by the column matrix
$[c_{\sigma }]^t_{|\sigma |=j}$, where $"t"$ denotes the 
matrix transpose. Then $T_{j,0}=0$ for 
all sufficiently large $j$'s and we can define 
$T\in \cL (\oplus _{k\geq 0}\CC _k)$
by using \eqref{baza}. Finally, set $\Phi _2(x)=T$.

\subsection{Spectral Factorization}
We briefly review the spectral factorization of positive definite 
kernels on the set $\NN _0$ of nonnegative integers. 
For more details, see \cite{Co}.
Let $\cE$ be a Hilbert space and let $\cP _+(\cE)$ be the set 
of positive definite kernels on $\NN _0$ with values in $\cL (\cE )$.
The order on $\cP _+(\cE)$ is: $K_1\leq K_2$ if $K_2-K_1$
belongs to $\cP _+(\cE)$. Next consider a family 
$\cF =\{\cF _n\}_{n\geq 0}$ of Hilbert spaces and call 
{\sl lower triangular array}
a family $\Theta =\{\Theta _{k,j}\}_{k,j\geq 0}$ of operators
$\Theta _{k,j}\in \cL (\cE ,\cF_k)$ such that $\Theta _{k,j}=0$
for $k<j$ and each column
$c_j(\Theta )=[\Theta _{k,j}]_{k\geq 0}$, $j\geq 0$,  belongs
to $\cL (\cE ,\oplus _{k\geq j}\cF _k)$. 
Denote by $H^2(\cE, \cF)$ the set of all lower triangular 
arrays as above. A lower triangular array
is called {\sl outer} if the set  
$\{c_j(\Theta )\cE\mid j\geq k\}$ is total in $\oplus _{j\geq k}\cF _j$
for all $k\geq 0$. 
If $\Theta $ is an outer triangular array, then 
the formula
$$K_\Theta (k,j)=c_k(\Theta )^*c_j(\Theta )$$
defines an element of $\cP _+(\cE)$.
For the proof of the following result see \cite{Co}, Chapter~5.

\begin{theorem}\label{factorizare}
Let $K$ be an element of $\cP _+(\cE)$. Then there exists a family
$\cF =\{\cF _n\}_{n\geq 0}$ of Hilbert spaces
and an outer triangular array $\Theta \in H^2(\cE ,\cF )$ such that 

\smallskip
$(a)$ $K_\Theta \leq K$.

\smallskip
$(b)$ For any other family 
$\cF '=\{\cF '_n\}_{n\geq 0}$ of Hilbert spaces
and any outer triangular array $\Theta '\in H^2(\cE ,\cF ')$ such that 
$K_{\Theta '}\leq K$, we have $K_{\Theta '}\leq K_{\Theta }$. 

\smallskip
$(c)$ $\Theta $ is uniquely determined by $(a)$ and $(b)$ up to a left 
unitary diagonal factor.
\end{theorem}

\subsection{Orthogonal Polynomials}
Let $\cP _{2N}$ be the algebra of polynomials
in $2N$ non-commuting
indeterminates $X_1$,$\ldots $,$X_N$, $X_{N+1}$,$\ldots $,$X_{2N}$ 
with complex coefficients.
An involution $\cJ$ can 
be introduced on $\cP _{2N}$ as follows:
$$\cJ (X_k)=X_{N+k},\quad k=1,\ldots ,N,$$
$$\cJ (X_l)=X_{l-N},\quad l=N+1,\ldots ,2N;$$
on monomials, 
$$\cJ (X_{i_1}\ldots X_{i_k})=\cJ (X_{i_k})\ldots \cJ (X_{i_1}),$$
and finally, if $Q=\sum _{\sigma \in \FF _{2N}^+}c_{\sigma }X_{\sigma }$,
then $\cJ (Q)=
\sum _{\sigma \in \FF _{2N}^+}\overline{c}_{\sigma }
\cJ (X_{\sigma })$.
Thus, $\cP _{2N}$ is a unital, associative, $*$-algebra over $\CC $
and we notice that $\cP _N$ is a subalgebra of $\cP _{2N}$.

We say that $\cA \subset \cP _{2N}$
is $\cJ $-symmetric if $P\in \cA $ implies 
$c\cJ (P)\in \cA $ for some $c\in \CC -\{0\}$. 
We construct an associative algebra
$\cT _N(\cA )$ as the quotient of $\cP _{2N}$
by the two-sided ideal $\cE (\cA )$ generated by 
$\cA $. We notice that $\cT _N(\emptyset )=
\cP _{2N}$.  We let $\pi =\pi _{\cA }:\cP _{2N}\rightarrow 
\cT _N(\cA )$ be the quotient map and since
$\cA $ is $\cJ $-symmetric,
\begin{equation}\label{invo}
\cJ _{\cA }(\pi  (P))=\pi (\cJ (P))
\end{equation}
gives an involution on $\cT _N(\cA )$.
We will be interested in linear functionals $\phi $
on $\cT _N(\cA )$ with the property that 
$\phi (\cJ _{\cA }(\pi (P))\pi (P))\geq 0$
for all $P\in \cP _N$ and we will say that $\phi $  
is a 
{\em positive functional} on $\cT _N(\cA )$.
We notice that 
$\phi (\cJ _{\cA }(\pi (P)))=
\overline{\phi (\pi (P))}$ for $P\in \cP _N$
and 
$$|\phi (\cJ _{\cA }(\pi (P_1))\pi (P_2))|^2\leq 
\phi (\cJ _{\cA }(\pi (P_1))\pi (P_1))
\phi (\cJ _{\cA }(\pi (P_2))\pi (P_2))
$$
for $P_1,P_2\in \cP _N$.

We now consider the GNS construction associated to
$\phi $. Thus, we define 
on $\pi (\cP _N)$,
\begin{equation}\label{tilde}
\langle \pi (P_1),\pi (P_2)\rangle _{\phi }=
\phi (\cJ _{\cA }(\pi (P_2)\pi (P_1)),
\end{equation}
and factor out the subspace
$\cN _{\phi }=\{\pi (P)\mid P\in \cP _N,
\langle \pi (P),\pi (P) \rangle _{\phi }=0\}$.
Completing this quotient with respect
to the norm induced by \eqref{tilde} we obtain a Hilbert space 
$\cH _{\phi }$.

From now on we will assume that 
$\phi $ 
is strictly positive, that is, 
$\phi (\cJ _{\cA }(\pi(P))\pi (P))>0$
for all $P\in \cP _N-\cE (\cA)$, so that  $\cN _{\phi }=\{0\}$ and 
$\pi (\cP _N)$ can be viewed as a subspace of $\cH _{\phi }$.
The {\em index set of $\cA$}, $G\subset \FF _N^+$, is  
chosen as follows:
let $\emptyset \in G$; if $\alpha \in G$, choose the next
element in $G$ to be the least $\beta \in \FF _N^+$ such that the set 
of elements $\pi (X_{\alpha '})$, $\alpha '\preceq \alpha $, 
and $\pi (X_{\beta })$ is linearly independent.
Define $F_{\alpha }=\pi (X_{\alpha })$ for $\alpha \in G$ and set
$\cB =\{F_{\alpha }\}_{\alpha \in G}$.
Let $G_n=\{\alpha \in G\mid |\alpha |=n\}$, then $G_0=
\{\emptyset \}$ and $\{G_n\}_{n\geq 0}$ is a partition of $G$.

Since $\phi $ is strictly positive it follows that 
$\cB $ is a linearly independent family in 
$\cH _{\phi }$ and 
the 
Gram-Schmidt procedure gives a family $\{\varphi _{\alpha  }\}
_{\alpha \in G}$ of elements in 
$\pi (\cP _N)\subset \cT _N(\cA )$ such that 

\begin{equation}\label{bond1}
\varphi _{\alpha  }=
\sum _{\beta \preceq \alpha }a_{\alpha ,\beta }F_{\beta },
\quad a_{\alpha ,\alpha }>0;
\end{equation}
\begin{equation}\label{bond2}
\langle \varphi _{\alpha }, \varphi _{\beta }\rangle _{\phi }=
\delta _{\alpha ,\beta },
\quad \alpha ,\beta \in G.
\end{equation}

\noindent
The elements $\varphi _{\alpha  }$, $\alpha \in G$, will be called
the {\em orthogonal polynomials} associated to $\phi $.
Typically, the theory of orthogonal polynomials deals
with the study of algebraic and asymptotic properties of
the orthogonal polynomials 
associated to strictly positive
functionals on $\cT _N(\cA )$.
We also notice that the use of the Gram-Schmidt process depends
on the order that we have chosen on 
$\FF _N^+$. A different order would yeald a different family of 
orthogonal polynomials. Due to the natural grading on 
$\FF _N^+$ it is possible to develop a base free approach to 
orthogonal polynomials. In the case of orthogonal polynomials
on several commuting variables this is presented in \cite{DX}. 
However, in this paper we stick to the lexicographic order on  
$\FF _N^+$ (and on the index set $G$).

An explicit formula for the orthogonal polynomials  
can be obtained in the same manner as in 
the classical (one scalar variable) case. Define
\begin{equation}\label{stilde}
s_{\alpha ,\beta }=
\phi (\cJ _{\cA }(F_{\alpha })F_{\beta })=
\langle F_{\beta },F_{\alpha }\rangle _{\phi }, 
\quad \alpha ,\beta \in G,
\end{equation}
and 
\begin{equation}\label{Deh}
D_{\alpha }=\det\left[s_{\alpha ',\beta '}\right]_
{\alpha ',\beta '\preceq \alpha }>0, \quad \alpha \in G.
\end{equation}  
We notice that $\phi $ is a positive functional on 
$\cT _N(\cA )$ if and only if $K_{\phi }(\alpha ,\beta )=
s_{\alpha ,\beta }$, $\alpha ,\beta \in G$, 
is a positive definite kernel on $G$.
While the kernel $K_{\phi }$ characterizes the positivity of the functional
$\phi $ and contains the basic information for the construction of the 
orthogonal polynomials, in general it does not determine $\phi $ uniquely.
We will occasionally say that the orthogonal polynomials are associated to the 
kernel $K_{\phi }$ rather then $\phi $ itself.
One typical situation when $K_{\phi }$ determines $\phi $ is when  
$\{\cJ(X_k)-X_k\mid k+1,\ldots ,N\}\subset \cA $; another example 
is provided by the Wick polynomials,
$$X_i\cJ (X_j)-\delta _{ij}+\sum _{k,l=1}^NT_{ij}^{kl}\cJ (X_l)X_k,
\quad i,j=1,\ldots ,N,$$
where $T_{ij}^{kl}$ are complex numbers and $\delta _{ij}$
is the Kronecker symbol (see \cite{JSW}). 

From now on $\tau -1$ denotes the predecessor of $\tau $ with respect to the 
lexicographic order on $\FF _N^+$, while $\sigma +1$
denotes the successor of $\sigma $. 
It is showed in \cite{Co} that 
$\varphi _{\emptyset }=s_{\emptyset ,\emptyset }^{-1/2}$ and for 
$\emptyset \prec \alpha $,
\begin{equation}\label{crop1}
\varphi _{\alpha }=\frac{1}{\sqrt{D_{\alpha -1}D_{\alpha }}}
{\det \left[\begin{array}{c}
\left[s_{\alpha ',\beta '}\right]_
{\alpha '\prec \alpha ;
\beta '\preceq \alpha } \\
  \\
\begin{array}{ccc}
F_{\emptyset } & \ldots & F_{\alpha }
\end{array}
\end{array}
\right]},
\end{equation}
with an appropriate interpretation of the determinant.
In most of the cases, the formula \eqref{crop1}
is not very useful for the actual computation of the orthogonal
polynomials or for their study. Instead there are used recurence 
formulae. We discuss several examples in the next sections.

\section{The case $\cA =\emptyset $, $N=1$}
It turns out that this is, in fact, the most general situation. 
For this reason we treat this case separately.
For $\cA =\emptyset $ and  $N=1$, the index set is $G=\NN _0$ and a linear 
functional on $\cP _2$ is positive if and only if 
$K_{\phi }(n,m)=\phi (\cJ (X_1^n)X_1^m)$, $n,m\in \NN _0$, is a positive
definite kernel on $\NN _0$. 
However, $K_{\phi }$ does not completely determine $\phi $.
Thus, there is no way to deduce $\phi (X_1\cJ (X_1))$ from $K_{\phi }$ in 
general. Still, we notice that there is no other restriction on $K_{\phi }$, 
in the sense that given a positive definite kernel $K$ on $\NN _0$, there
exist positive functionals $\phi $ on $\cP _2$ such that $K_{\phi }=K$.
This is done simply by the linearization of any function 
$\phi _0:\FF _{2}^+\rightarrow \CC $ such that 
$\phi _0(\cJ (X_{\alpha })X_{\beta })=K(\alpha ,\beta )$
for $\alpha ,\beta \in \FF _2^+$.  

Let $\{\gamma _{k,j}\}_{0\leq k<j}$ be the parameters 
associated to $K_{\phi }$ by 
\cite{Co}, Theorem~1.5.3. Assuming that $\phi $ is strictly positive 
means $|\gamma _{k,j}|<1$ for all $0\leq k<j$. Define 
$d_{k,j}=(1-|\gamma _{k,j}|^2)^{1/2}$.  
An explicit connection between $K_{\phi }$ and $\{\gamma _{k,j}\}_{0\leq k<j}$
is given by formula (1.4.6) in \cite{Co},
\begin{equation}\label{baza}
s_{k,j}=s^{1/2}_{k,k}\left[
\begin{array}{ccc}
1 & 0 & \ldots 
\end{array}
\right]
U_{k,j}
\left[
\begin{array}{c}
1 \\
 0 \\
 \vdots 
\end{array}
\right]s^{1/2}_{j,j},
\end{equation}
where $U_{k,j}$ is defined by (1.5.5) in \cite{Co}. 
For a better understanding of this formula it could be useful
to write it explicitely for a few particular
indices and to answer a related combinatorial question.
Thus, 
$$s_{01}=s^{1/2}_{00}\left[
\begin{array}{cc}
1 & 0   
\end{array}
\right]
\left[
\begin{array}{cc}
\gamma _{01} & d_{01} \\
d_{01} & -\overline{\gamma }_{01} 
\end{array}
\right]
\left[
\begin{array}{c}
1 \\
 0 
\end{array}
\right]s^{1/2}_{11}=s^{1/2}_{00}\gamma _{01}s^{1/2}_{11};
$$
$$\begin{array}{rcl}
s_{02}\!\!\!\!& =&\!\!\!\!s^{1/2}_{00}\left[
\begin{array}{ccc}
1 & 0 & 0  
\end{array}
\right]
\left[
\begin{array}{ccc}
\gamma _{01} & d_{01} & 0 \\
d_{01} & -\overline{\gamma }_{01} & 0 \\
 0 & 0 & 1
\end{array}
\right]
\left[
\begin{array}{ccc}
 1 & 0 & 0 \\
 0 & \gamma _{02} & d_{02} \\
 0 & d_{02} & -\overline{\gamma }_{02} 
\end{array}
\right]
\left[
\begin{array}{ccc}
\gamma _{12} & d_{12} & 0 \\
d_{12} & -\overline{\gamma }_{12} & 0 \\
 0 & 0 & 1
\end{array}
\right]
\left[
\begin{array}{c}
1 \\
 0 \\
 0 
\end{array}
\right]s^{1/2}_{22} \\
 & & \\
& =&\!\!\!\!s^{1/2}_{00}\left(
\gamma _{01}\gamma _{12}+d_{01}\gamma _{02}d_{12}\right)s^{1/2}_{22};
\end{array}
$$
$$\begin{array}{rcl}
s_{03}\!\!\!\!&=&\!\!\!\!s^{1/2}_{00}\left[
\begin{array}{cccc}
1 & 0 & 0 & 0 
\end{array}
\right]
\left[
\begin{array}{cccc}
\gamma _{01} & d_{01} & 0 & 0 \\
d_{01} & -\overline{\gamma }_{01} & 0 & 0 \\
 0 & 0 & 1 & 0 \\
 0 & 0 & 0 & 1
\end{array}
\right]
\left[
\begin{array}{cccc}
 1 & 0 & 0 & 0 \\
 0 & \gamma _{02} & d_{02} & 0 \\
 0 & d_{02} & -\overline{\gamma }_{02} & 0 \\
  0 & 0 & 0 & 1
\end{array}
\right]
\left[
\begin{array}{cccc}
 1 & 0 & 0 & 0 \\
 0 & 1 & 0 & 0 \\
 0 & 0 & \gamma _{03} & d_{03} \\
 0 & 0 & d_{03} & -\overline{\gamma }_{03} \\
\end{array}
\right] \\
 & & \\
 & & \,\,\,\,\,\,\,\, \quad \quad \quad \quad \quad \quad \left[ 
\begin{array}{cccc}
\gamma _{12} & d_{12} & 0 & 0 \\
d_{12} & -\overline{\gamma }_{12} & 0 & 0 \\
 0 & 0 & 1 & 0 \\
 0 & 0 & 0 & 1
\end{array}
\right]
\left[
\begin{array}{cccc}
 1 & 0 & 0 & 0 \\
 0 & \gamma _{13} & d_{13} & 0 \\
 0 & d_{13} & -\overline{\gamma }_{13} & 0 \\
  0 & 0 & 0 & 1
\end{array}
\right] \\
 & & \\
 & & \,\,\,\,\,\,\,\, \quad \quad \quad \quad \quad \quad \left[
\begin{array}{cccc}
\gamma _{23} & d_{23} & 0 & 0 \\
d_{23} & -\overline{\gamma }_{23} & 0 & 0 \\
 0 & 0 & 1 & 0 \\
 0 & 0 & 0 & 1
\end{array}
\right]
\left[
\begin{array}{c}
1 \\
 0 \\
 0 \\
 0
\end{array}
\right]s^{1/2}_{33} \\
 & & \\
 &=&\!\!\!\!s^{1/2}_{00}\left(\gamma _{01}\gamma _{12}\gamma _{23}+
\gamma _{01}d_{12}\gamma _{13}d_{23}+
d_{01}\gamma _{02}d_{12}\gamma _{23} 
-d_{01}\gamma _{02}\overline{\gamma }_{12}\gamma _{13}d_{23}+
d_{01}d_{02}\gamma _{03}d_{13}d_{23}\right)s^{1/2}_{33}.
\end{array}
$$
A natural combinatorial question would be to calculate the number 
$N(s_{k,j})$ of additive terms in the expression of $s_{k,j}$.
Thus, for $k\geq 0$,
$$N(s_{01})=N({s_{k,k+1}})=1,$$
$$N(s_{02})=N({s_{k,k+2}})=2,$$
$$N(s_{03})=N({s_{k,k+3}})=5.$$
The general formula is given by the following result.
\begin{theorem}\label{catalan}
$N(s_{k,k+l})$ is given by the Catalan number $\displaystyle\frac{1}{l+1}
\left(\begin{array}{c}
2l \\
l 
\end{array}
\right)$.
\end{theorem}
\begin{proof}
The first step of the proof considers the realization of $s_{k,j}$
through a time varying transmission line (or lattice) (see \cite{Co}, 
Chapter~4,
for more details). For illustration we consider the case of $s_{03}$
in Figure ~1.

\begin{figure}[h]
\setlength{\unitlength}{2700sp}%
\begingroup\makeatletter\ifx\SetFigFont\undefined%
\gdef\SetFigFont#1#2#3#4#5{%
  \reset@font\fontsize{#1}{#2pt}%
  \fontfamily{#3}\fontseries{#4}\fontshape{#5}%
  \selectfont}%
\fi\endgroup%
\begin{picture}(6324,2874)(289,-2323)
{ \thinlines
%\put(601,-361){\circle{300}}
}%
%{ \put(2101,-361){\circle{300}}
%}%
%{ \put(3901,-361){\circle{300}}
%}%
%{ \put(5701,-361){\circle{300}}
%}%
{ \put(301,-361){\line( 1, 0){900}}
}%
%{ \put(601,-211){\line( 0,-1){300}}
%}%
{ \put(1201,-361){\vector( 1,-1){600}}
}%
{ \put(1201,-961){\vector( 1, 1){600}}
}%
%{ \put(451,239){\framebox(300,300){ L_{44}}}
%}%
%{ \put(601,239){\vector( 0,-1){450}}
%}%
%{ \put(1051,-361){\line( 0, 1){300}}
%}%
%{ \put(1051,239){\vector( 0, 1){300}}
%}%
%{ \put(901,-61){\framebox(300,300){L^*_{44}}}
%}%
{ \put(500,10){A}}
{ \put(6000,10){B}}
{ \put(1126,-1111){\framebox(750,900){}}
}%
{ \put(901,-961){\line( 1, 0){1200}}
}%
{ \put(1201,-361){\line( 1, 0){1200}}
}%
{ \put(2926,-1111){\framebox(750,900){}}
}%
{ \put(2026,-1711){\framebox(750,900){}}
}%
{ \put(3826,-1711){\framebox(750,900){}}
}%
{ \put(4726,-1111){\framebox(750,900){}}
}%
{ \put(2926,-2311){\framebox(750,900){}}
}%
{ \put(2101,-961){\line( 1, 0){1200}}
}%
{ \put(3301,-961){\line( 1, 0){1200}}
}%
{ \put(4501,-961){\line( 1, 0){1200}}
}%
{ \put(2401,-361){\line( 1, 0){1200}}
}%
{ \put(3601,-361){\line( 1, 0){1200}}
}%
{ \put(4801,-361){\line( 1, 0){1200}}
}%
{ \put(1201,-361){\vector( 1, 0){600}}
}%
{ \put(1201,-961){\vector( 1, 0){600}}
}%
{ \put(3001,-361){\vector( 1, 0){600}}
}%
{ \put(3001,-961){\vector( 1, 0){600}}
}%
{ \put(2101,-961){\vector( 1, 0){600}}
}%
{ \put(2101,-1561){\vector( 1, 0){600}}
}%
{ \put(3001,-1561){\vector( 1, 0){600}}
}%
{ \put(3001,-2161){\vector( 1, 0){600}}
}%
{ \put(3901,-961){\vector( 1, 0){600}}
}%
{ \put(3901,-1561){\vector( 1, 0){600}}
}%
{ \put(4801,-361){\vector( 1, 0){600}}
}%
{ \put(4801,-961){\vector( 1, 0){600}}
}%
{ \put(3001,-961){\vector( 1, 1){600}}
}%
{ \put(4801,-961){\vector( 1, 1){600}}
}%
{ \put(3001,-361){\vector( 1,-1){600}}
}%
{ \put(4801,-361){\vector( 1,-1){600}}
}%
{ \put(2101,-1561){\vector( 1, 1){600}}
}%
{ \put(3901,-1561){\vector( 1, 1){600}}
}%
{ \put(3001,-2161){\vector( 1, 1){600}}
}%
{ \put(2176,-961){\vector( 1,-1){600}}
}%
{ \put(3001,-1561){\vector( 1,-1){600}}
}%
{ \put(3901,-961){\vector( 1,-1){600}}
}%
{ \put(1801,-1561){\line( 1, 0){3000}}
}%
{ \put(2701,-2161){\line( 1, 0){1200}}
}%
%{ \put(2101,-511){\line( 0, 1){300}}
%}%
%{ \put(3901,-511){\line( 0, 1){300}}
%}%
%{ \put(5701,-511){\line( 0, 1){300}}
%}%
%{ \put(1951,239){\framebox(300,300){L_{33}}}
%}%
%{ \put(3751,239){\framebox(300,300){L_{22}}}
%}%
%{ \put(5551,239){\framebox(300,300){L_{11}}}
%}%
%{ \put(2101,239){\vector( 0,-1){450}}
%}%
%{ \put(3901,239){\vector( 0,-1){450}}
%}%
%{ \put(5701,239){\vector( 0,-1){450}}
%}%
%{ \put(2401,-61){\framebox(300,300){L^*_{33}}}
%}%
%{ \put(4201,-61){\framebox(300,300){L^*_{22}}}
%}%
%{ \put(2551,-361){\line( 0, 1){300}}
%}%
%{ \put(4351,-361){\line( 0, 1){300}}
%}%
%{ \put(2551,239){\vector( 0, 1){300}}
%}%
%{ \put(4351,239){\vector( 0, 1){300}}
%}%
%{ \put(6001,-511){\framebox(300,300){L^*_{11}}}
%}%
%{ \put(6376,-361){\vector( 1, 0){225}}
%}%
\end{picture}

\caption{\mbox{ Lattice representation for $s_{03}$}}
\end{figure}

Each box in Figure~1 represents the action of the unitary matrix
$$\left[\begin{array}{cc}
\gamma _{k,j} & d_{k,j} \\
d_{k,j} & -\overline{\gamma }_{k,j}
\end{array}\right]
$$
and we see that the number of additive terms in the formula of 
$s_{03}$ is given by
the number of paths from $A$ to $B$ in Figure~1. In it clear that 
to each 
path from $A$ to $B$ in Figure ~1 it corresponds a ``good'' path from 
$C$ to $D$ in Figure ~2, that is, a path that never steps 
below the diagonal and goes only to the right or 
downward. 

\begin{figure}[h]
\setlength{\unitlength}{2700sp}%
\begingroup\makeatletter\ifx\SetFigFont\undefined%
\gdef\SetFigFont#1#2#3#4#5{%
  \reset@font\fontsize{#1}{#2pt}%
  \fontfamily{#3}\fontseries{#4}\fontshape{#5}%
  \selectfont}%
\fi\endgroup%
\begin{picture}(1800,1733)(76,-1111)
{\thinlines
\put(751,539){\circle{150}}
}%
{ \put(1201,539){\circle{150}}
}%
{\put(-100,539){{C}}
}%
{\put(1900,-811){{D}}
}%
{ \put(1651,539){\circle{150}}
}%
{ \put(301, 89){\circle{150}}
}%
{ \put(751, 89){\circle{150}}
}%
{ \put(1201, 89){\circle{150}}
}%
{ \put(1651, 89){\circle{150}}
}%
{ \put(301,-361){\circle{150}}
}%
{ \put(751,-361){\circle{150}}
}%
{ \put(1201,-361){\circle{150}}
}%
{ \put(1651,-361){\circle{150}}
}%
{ \put(301,-811){\circle{150}}
}%
{ \put(751,-811){\circle{150}}
}%
{ \put(1201,-811){\circle{150}}
}%
{ \put(1651,-811){\circle{150}}
}%
{ \put(301,539){\circle{150}}
}%
{ \put(301,539){\line( 1, 0){450}}
\put(751,539){\line( 0,-1){450}}
\put(751, 89){\line( 1, 0){900}}
\put(1651, 89){\line( 0,-1){900}}
}%
\end{picture}

\caption{\mbox{ A ``good'' path from $C$ to $D$}}
\end{figure}

More precisely, each box in Figure~1 corresponds to a point 
strictly above the diagonal 
in Figure~2. Once this one-to-one correspondence is established, 
we can use the 
well-known fact that the number of ``good'' paths like the one in Figure~2 
is given 
exactly by the Catalan numbers.
\end{proof}

Returning to orthogonal polynomials, we notice that they
obey the following recurrence 
relations (see \cite{CJ}, formulae (3.10) and (3.11)):
\begin{equation}\label{lazero}
\varphi _0(X,l)=\varphi _0^{\sharp }(X,l)=s_{l,l}^{-1/2}, \quad l\in \NN _0,
\end{equation}
and for $n\geq 1$, $l\in \NN _0$,
\begin{equation}\label{primarelatie}
\varphi _n(X,l)=\frac{1}{d_{l,n+l}}
\left( X\varphi _{n-1}(X,l+1)-
\gamma _{l,n+l}\varphi ^{\sharp }_{n-1}(X,l)\right),
\end{equation}
\begin{equation}\label{adouarelatie}
\varphi ^{\sharp}_n(X,l)=\frac{1}{d_{l,n+l}}
\left(-\overline{\gamma }_{l,n+l}X\varphi _{n-1}(X,l+1)+
\varphi ^{\sharp }_{n-1}(X,l)\right),
\end{equation}
where $\varphi _n(X)=\varphi _n(X,0)$ and 
$\varphi ^{\sharp }_n(X)=\varphi ^{\sharp }_n(X,0)$.
While it is clear how to recover the coefficients $\gamma _{k,j}$
from $K_{\phi }$, it appears to be also useful to recover these
parameters from the orthogonal polynomials.
It follows from the proof of Theorem ~3.2 in \cite{CJ} that 
$\{\varphi _n(X,l)\}_{n\geq 0}$ is the family of orthogonal polynomials
associated to the kernel $K^l_{\phi }(\alpha ,\beta )=
s_{\alpha +l,\beta +l}$,
$\alpha , \beta \in \NN _0$. Let $k_n^l$ be the leading coefficient
of $\varphi _n(X,l)$. We obtain the following formula for the parameters
$\gamma _{k,j}$.
\begin{theorem}\label{nouaparam}
For $l\in \NN _0$ and $n\geq 1$, 
$$\gamma _{l,n+l}=-s^{1/2}_{00}s^{-1/2}_{l+1,l+1}
\varphi _n(0,l)\displaystyle\frac{k^{l+1}_1\ldots k^{l+1}_{n-1}}
{k^l_1\ldots k^l_n}.
$$
\end{theorem}
\begin{proof}
We deduce
from \eqref{primarelatie} that
$$\varphi _n(0,l)=-\frac{\gamma _{l,n+l}}{d_{l,n+l}}
\varphi ^{\sharp}_{n-1}(0,l),
$$
and from formula \eqref{adouarelatie} 
we deduce 
$$\varphi ^{\sharp }_n(0,l)=
\frac{1}{d_{l,n+l}}\varphi ^{\sharp }_{n-1}(0,l)=\ldots 
=s_{0,0}^{-1/2}\prod _{p=1}^n\frac{1}{d_{l,p+l}},
$$
hence
$$\varphi _n(0,l)=-s_{0,0}^{-1/2}\gamma _{l,n+l}
\prod _{p=1}^n\frac{1}{d_{l,p+l}}.
$$
Using Theorem ~1.5.10 in \cite{Co}, we deduce that 
$$\prod _{p=1}^nd^2_{l,p+l}=s^{-1}_{l,l}
\displaystyle\frac{D_{l,l+n}}{D_{l+1,l+n}}
$$
so that, 
\begin{equation}\label{gama}
\gamma _{l,n+l}=-s^{1/2}_{00}s^{-1/2}_{l,l}
\varphi _n(0,l)\sqrt{\frac{D_{l,l+n}}{D_{l+1,l+n}}}.
\end{equation}
On the other hand, \eqref{primarelatie} gives that
$$k^l_n=\prod _{p=1}^{n-1}\frac{1}{d_{l+p,l+n}},
$$
and using once again  Theorem ~1.5.10 in \cite{Co}, we deduce
$$k^l_n=\sqrt{\frac{D_{l,l+n-1}}{D_{l,l+n}}}.
$$
This implies that
$$k^l_1\ldots k^l_n=\displaystyle\frac{s^{1/2}_{l,l}}{\sqrt{D_{l,l+n}}},
$$
and this can be used in \eqref{gama} in order to conclude the 
proof.
\end{proof}

We now develop an analogue of \eqref{unu}
and \eqref{doi}. The formulae
\eqref{primarelatie} and \eqref{adouarelatie} suggest that it is 
more convenient to work in a larger algebra. Thus, we consider the 
set $\cR _{1}$ of lower triangular arrays 
$a=[a_{k,j}]_{k,j\geq 0}$ with complex entries.
No boundedness assumption is made on these arrays. The addition in 
$\cR _{1}$ is defined by entry-wise addition and the multiplication is 
defined as follows: for $a=[a_{k,j}]_{k\geq j}$, 
$b=[b_{k,j}]_{k,j\geq 0}$ two elements of $\cR _{1}$,
$$(ab)_{k,j}=\sum _{l\geq 0}a_{k,l}b_{l,j}$$
(the sum is finite since both $a$ and $b$ are lower triangular). Thus,
$\cR _{1}$ becomes an associative, unital algebra. 

Next we associate the element $\Phi _n$ of $\cR _{1}$  
to the polynomials 
$\varphi _n(X,l)=\sum _{k=0}^na_{n,k}^lX^k$, $n,l\geq 0$, by the formula
\begin{equation}\label{fiunu}
(\Phi _n)_{k,j}=\left\{\begin{array}{lcl}
a_{n,k-j}^j & \quad & k\geq j \\
0    & \quad & k<j;
\end{array}\right.
\end{equation}
similarly, the element $\Phi ^{\sharp }_n$ of $\cR _{1}$ is associated 
to the family of polynomials 
$\varphi ^{\sharp }_n(X,l)=\sum _{k=0}^nb_{n,k}^lX^k$, $n,l\geq 0$, 
by the formula
\begin{equation}\label{fidoi}
(\Phi ^{\sharp }_n)_{k,j}=\left\{\begin{array}{lcl}
b_{n,k-j}^j &\quad  & k\geq j \\
0 & \quad  & k<j.
\end{array}\right.
\end{equation}

We notice that since $K_{\phi }$ is a scalar-valued kernel, the Hilbert
spaces $\cF _n$, $n\geq 0$, given by Theorem ~\ref{factorizare} are at most
one-dimensional (see \cite{Co}, Section ~5.1 for details). This implies that 
we can uniquely determine the spectral factor $\Theta _{\phi }$ of $K_{\phi }$ 
by the requirement that $(\Theta _{\phi })_{n,n}\geq 0$ for all $n\geq 0$.
Also, $\Theta _{\phi }\in \cR _{1}$. From now on we assume that 
$(\Theta _{\phi })_{n,n}>0$ for all $n\geq 0$ 
and we say that in this case $\phi $ (or $K_{\phi }$) belongs to the 
Szeg\"o class. By formula (5.1.5) in \cite{Co}, it follows that $\phi $ 
belongs to the Szeg\"o class if and only if 
\begin{equation}\label{prod}
s_{k,k}\prod _{n>k}d_{k,n}>0
\end{equation}
for all $k\geq 0$. This implies, in particular, that 
$\Phi ^{\sharp }_n$ is invertible in $\cR _{1}$ for all $n\geq 0$.
Finally, we say that a sequence $\{a_n\}\subset  \cR _{1}$
converges to $a\in \cR _{1}$ if $\{(a_n)_{k,j}\}$ 
converges to $a_{k,j}$ for all $k,j\geq 0$ (and we write $a_n\rightarrow a$).
We now obtain the following generalization of \eqref{unu} and 
\eqref{doi}.

\begin{theorem}\label{convergenta}
Let $\phi $ belong to the Szeg\"o class. Then 

\begin{equation}\label{doiunu}
\Phi _n\rightarrow 0
\end{equation}
and 
\begin{equation}\label{doidoi}
(\Phi _n^{\sharp })^{-1}\rightarrow \Theta _{\phi }.
\end{equation}
\end{theorem}

\begin{proof}
First we show \eqref{doiunu}. It is convenient to consider 
the natural derivation on $\cP _1$; for $P=\sum _{k=0}^na_kX^k\in \cP _1$, 
$$P^{(1)}=\sum _{k=1}^nka_kX^{k-1},$$
and then, for $k\geq 1$, 
$$P^{(k)}=(P^{(k-1)})^{(1)}.$$
We see that \eqref{doiunu} is equivalent to 
$$\varphi ^{(k)}_n(0,l)\rightarrow 0$$
for each fixed $k,l\geq 0$. We claim that
\begin{equation}\label{doitrei}
\sum _{n\geq 0}|\varphi ^{(k)}_n(0,l)|^2<\infty 
\end{equation}
and for each $l,k\geq 0$, 
\begin{equation}\label{doipatru}
\lim _{n\rightarrow \infty }(\varphi ^{\sharp }_n)^{(k)}(0,l)
\quad \mbox{exists and is finite.}
\end{equation}

We prove these statements by induction on $k\geq 0$. For $k=0$ we use
the formula
$$\varphi _n(0,l)=-s_{0,0}^{-1/2}\gamma _{l,n+l}
\prod _{p=1}^n\frac{1}{d_{l,p+l}}.
$$
obtained in the proof of Theorem ~\ref{nouaparam}
in order to deduce that 
$$\sum _{n\geq 0}|\varphi ^{(k)}_n(0,l)|^2=
s_{0,0}^{-1}\sum _{n\geq 0}|\gamma _{l,n+l}|^2
\prod _{p=1}^n\frac{1}{d^2_{l,p+l}}.
$$
Since $\phi $ belongs to the Szeg\"o class, we have that
$$g_l=s_{l,l}\prod _{n>l}d_{l,n}>0,$$
hence $\prod _{p=1}^n\frac{1}{d^2_{l,p+l}}\leq c_l$ for some 
$c_l>0$ and all $n\geq 0$. Also, for all $n\geq 0$, 
$$\sum _{n\geq 0}|\gamma _{l,n+l}|^2<\infty .$$
In particular, $\gamma _{l,n+l}\rightarrow 0$ as $n\rightarrow \infty $;
all of these give \eqref{doitrei} and \eqref{doipatru} for 
$k=0$ and all $l\geq 0$.

We now proceed to prove the general case. From \eqref{primarelatie} and 
\eqref{adouarelatie} we also deduce 
\begin{equation}\label{treiunu}
\varphi ^{(k)}_n(0,l)=\frac{1}{d_{l,n+l}}
\left(k\varphi ^{(k-1)}_{n-1}(0,l+1)-
\gamma _{l,n+l}(\varphi ^{\sharp }_{n-1})^{(k)}(0,l)\right),
\end{equation}
\begin{equation}\label{treidoi}
(\varphi ^{\sharp }_n)^{(k)}(0,l)=\frac{1}{d_{l,n+l}}
\left(-k\overline{\gamma }_{l,n+l}\varphi ^{(k-1)}_{n-1}(0,l+1)+
(\varphi ^{\sharp }_{n-1})^{(k)}(0,l)\right)
\end{equation} 
for $k\geq 1$.

Since $k\geq 1$, $(\varphi ^{\sharp }_0)^{(k)}(0,l)=0$, and we deduce 
from \eqref{treidoi} that 
$$(\varphi ^{\sharp }_n)^{(k)}(0,l)=
-k\left(\prod _{p=1}^n\frac{1}{d_{l,p+l}}\right)\sum _{j=1}^n
\overline{\gamma }_{l,j+l}\left(\prod _{q=1}^jd_{l,q+l}\right)
\varphi ^{(k-1)}_{j-1}(0,l+1).
$$
By Schwarz inequality, 
$$\begin{array}{l}
\sum _{j\geq 1}
|\overline{\gamma }_{l,j+l}\left
(\prod _{q=1}^jd_{l,q+l}\right)\varphi ^{(k-1)}_{j-1}(0,l+1)| \\
 \\
\quad \leq \left(\sum _{j\geq 1}|\gamma _{l,j+l}|^2
\prod _{q=1}^jd^2_{l,q+l}\right)^{1/2}
\left(\sum _{j\geq 1}
|\varphi ^{(k-1)}_{j-1}(0,l+1)|^2\right)^{1/2}.
\end{array}
$$
Again, since $\phi $ belongs to the Szeg\"o class and 
$\prod _{q=1}^jd^2_{l,q+l}\leq 1$, we deduce that 
$$\begin{array}{l}
\sum _{j\geq 1}
|\overline{\gamma }_{l,j+l}\varphi ^{(k-1)}_{j-1}(0,l+1)
\prod _{q=1}^jd_{l,q+l}| \\
 \\
\quad \leq C\left(\sum _{j\geq 1}
|\varphi ^{(k-1)}_{j-1}(0,l+1)|^2\right)^{1/2}.
\end{array}
$$
This and the induction hypothesis give that the series
$$\sum _{j\geq 1}
\overline{\gamma }_{l,j+l}\left(\prod _{q=1}^jd_{l,q+l}\right)
\varphi ^{(k-1)}_{j-1}(0,l+1)
$$
converges absolutely and since 
$$
\lim _{n\rightarrow \infty }\prod _{p=1}^n\frac{1}{d_{l,p+l}}=
\frac{s_{l,l}}{g_l}<\infty ,
$$
we deduce that 
$\lim _{n\rightarrow \infty }(\varphi ^{\sharp }_n)^{(k)}(0,l)$
exists and is finite.

Using \eqref{treiunu},
$$\begin{array}{l}
\sum _{n\geq 1}
|\varphi ^{(k)}_{n}(0,l)|^2 \\
 \\
\,\,\,\leq k^2\sum _{n\geq 1}\frac{1}{d_{l,n+l}^2}
|\varphi ^{(k-1)}_{n-1}(0,l+1)|^2 \\
 \\
\,\,\,\,\,\, +2k\sum _{n\geq 1}\frac{1}{d_{l,n+l}^2}
|\varphi ^{(k-1)}_{n-1}(0,l+1)\gamma _{l,n+l}
(\varphi ^{\sharp }_{n-1})^{(k)}(0,l)| \\
 \\
\,\,\,\,\,\,\,\,\, +\sum _{n\geq 1}\frac{1}{d_{l,n+l}^2}
|\gamma _{l,n+l}|^2|(\varphi ^{\sharp }_{n-1})^{(k)}(0,l)|^2.
\end{array}
$$

Since for sufficiently large $n$, $d_{l,n+l}\geq C_l>0$
and $|(\varphi ^{\sharp }_{n-1})^{(k)}(0,l)|\leq C'_l$, another 
application of the Schwarz inequality, the fact that $\phi $ 
belongs to the Szeg\"o class, and the induction hypothesis give that 
$\sum _{n\geq 1}
|\varphi ^{(k)}_{n}(0,l)|^2<\infty $. In particular, 
$\varphi _n^{(k)}(0,l)\rightarrow 0$ as $n\rightarrow \infty $, 
concluding the proof of 
\eqref{doiunu}.

 A convenient proof of \eqref{doidoi} can be based on the so-called 
Toeplitz embedding, systematically used in \cite{FFGK}.
This approach would also explain the meaning of the elements 
$\Phi _n$, $\Phi _n^{\sharp }$ of $\cR _1$. Define, for $n\geq 1$, 
$$
(\Gamma _n)_{k,j}=\left\{\begin{array}{lcl}
\gamma _{k,j} &\quad  & j=k+n \\
0 & \quad  & \mbox{otherwise}.
\end{array}\right.
$$
Then let $A=\left[A_{j-k}\right]_{k,j\geq 0}$ be the positive definite 
Toeplitz kernel associated by Proposition~1.5.6 in \cite{Co} 
to $\{\Gamma _n\}_{n\geq 1}$ and  
$A_0=\left[s_{l,l}\right]$. By Proposition~1.6.10 (a) in \cite{Co},
$K_{\phi }$ is just a compression of the kernel $A$.
By Proposition~1.6.10 (b) and Theorem ~5.1.2 in \cite{Co}, 
the spectral factor $\Theta _{\phi }$ of $K_{\phi }$ is a corresponding 
compression of the spectral factor of $A$. The key point of the proof 
is the connection between $\Phi _n$ and the right orthogonal polynomials of 
$A$. For details on operator-valued orthogonal polynomials see
\cite{BC}. Thus, let $\{R_n\}_{n\geq 0}$ be the set of the right 
orthogonal polynomials of $A$, 
$$R_n(z)=\sum _{k=0}^nR_{n,k}z^k, \quad R_{n,n}\geq 0,$$
and define $R^{\sharp }_n=z^nR_n(1/\overline{z})^*=
\sum _{k=0}^nnR_{n,n-k}z^k$. Also, define 
$$\rho _n=\left[R^*_{n,n},\ldots ,R^*_{n,0}\right]^t$$
and let $\tilde \rho _n$ be obtained by canonical reshuffle of $\rho _n$
(\cite{Pa}, Chapter~7). Then, 
$$\Phi ^{\sharp }_n=\left[
\begin{array}{cc}
I_n & 0 \\
0 & \tilde \rho _n 
\end{array}
\right],
$$
where $I_n$ denotes the $n\times n$ identity matrix. This relation 
can be easily checked by using the characterization of 
$$\left[R_{n,0},\ldots ,R_{n,n}\right]^t$$
as the unique solution of 
$$\left[A_{j-k}\right]_{0\leq k,j\leq n}
\left[R_{n,0},\ldots ,R_{n,n}\right]^t=
\left[0,\ldots ,0,D_n\right]^t,
$$
where $D_n$ is a positive operator, and the orthogonality properties
of $\{\phi _n\}_{n\geq 0}$. Finally, an application of Theorem~4.37 in 
\cite{BC} concludes the proof.
\end{proof}

In order to provide generalizations of \eqref{trei} and \eqref{patru}
in this setting we consider first their geometrical interpretation.
Thus, by a result of Kolmogorov, $K_{\phi }$ is the covariance 
kernel of a stochastic process $\{f_n\}_{n\geq 0}\subset L^2(\mu )$
for some probability space $(X,\cM ,\mu )$. That is, 
$$K_{\phi }(m,n)=\int _Xf_n\overline{f}_md\mu .$$
We can suppose, without loss of generality, that $\{f_n\}_{n\geq 0}$
is total in $L^2(\mu )$ and for $p\leq q$ we introduce 
the subspaces $\cE _{p,q}$ 
given by the closure in $L^2(\mu )$ of the linear span of 
$\{f_k\}_{k=p}^q$. 

The operator angle between two spaces $\cE _1$ and $\cE _2$ of 
$L^2(\mu )$ is defined by 
$$B(\cE _1, \cE _2)=P_{\cE _1}P_{\cE _2}P_{\cE _1},$$
where $P_{\cE _1}$ is the orthogonal projection of 
$L^2(\mu )$ onto $\cE _1$. Also define 
$$\Delta (\cE _1, \cE _2)=I-B(\cE _1, \cE _2).$$

We associate to the process $\{f_n\}_{n\geq 0}$
a family of subspaces $\cH _{r,q}$ of $L^2(\mu )$
such that $\cH _{r,q}$ is the closure of the linear space generated
by $f_k$, $r\leq k\leq q$. 

The geometric interpretation of the limits 
\eqref{trei} and \eqref{patru}
is discussed in \cite{GS} and nonstationary extensions are
presented in \cite{Co}, Chapter~6. The interpretation of the second 
Szeg\"o limit theorem in \cite{Co} required a stochastic process 
indexed by the 
set of integers, 
which is not the case in our situation. So, we need a modification 
of that interpretation that fits into our setting. Thus, we consider 
first the scale of limits:
\begin{equation}\label{alfa}
s-\lim _{r\rightarrow \infty }\Delta (\cH _{0,n},\cH _{n+1,r})=
\Delta (\cH _{0,n},\cH _{n+1,\infty })
\end{equation}
for $n\geq 0$, and then we let $n\rightarrow \infty $ and deduce
\begin{equation}\label{beta}
s-\lim _{n\rightarrow \infty }\Delta (\cH _{0,n},\cH _{n+1,\infty })=
\Delta (\cH _{0,\infty },\cap _{n\geq 0}\cH _{n,\infty }),
\end{equation} 
where $s-\lim $ denotes the strong operatorial limit. 

We then deduce analogues of the Szeg\"o limit theorems by expressing 
these limits in terms of the determinants 
$D_{r,q}=\det \left[K_{\phi }(r',q')\right]_{r\leq r',q'\leq q}$, 
$r\leq q$.

\begin{theorem}\label{szego}
Let $\phi $ belong to the Szeg\"o class.
Then
\begin{equation}\label{3?}
\frac{D_{r,q}}{D_{r+1,q}}=s_{r,r}\det \Delta (\cH _{r,r},\cH _{r+1,q})=
\frac{s_{r,r}}{|\varphi ^{\sharp }_q(0,r)|^2},
\end{equation}
\begin{equation}\label{3??}
\lim _{q\rightarrow \infty }\frac{D_{r,q}}{D_{r+1,q}}=
s_{r,r}\det \Delta (\cH _{r,r},\cH _{r+1,\infty })=
|\Theta _{\phi }(r,r)|^2=s_{r,r}\prod _{j\geq 1}d^2_{r,r+j}.
\end{equation}
If we denote the above limit by $g_r$ and 
$$L=\lim _{n\rightarrow \infty }\prod _{0\leq k\leq n\leq j}
d_{k,j}^2>0,$$
then 
\begin{equation}\label{3???}
\lim _{n\rightarrow \infty }
\frac{D_{0,n}}{\prod _{l=0}^ng_l}
=\frac{1}{\det \Delta (\cH _{0,\infty }, \cap _{n\geq 0}\cH _{n,\infty })}
=\frac{1}{L}.
\end{equation}
\end{theorem}

\begin{proof}
The connection between the operator angles and determinants of 
type $D_{r,q}$ is given by the following formula
which is a consequence of Lemma ~6.4.1 in \cite{Co}: for 
$r\leq l\leq q,$
\begin{equation}\label{formula}
\det \Delta (\cH _{r,l},\cH _{l+1,q})=
\frac{D_{r,q}}{D_{r,l}D_{l+1,q}}.
\end{equation}
Then Theorem~1.5.10 in \cite{Co} allows the computation of $D_{r,q}$ in terms
of the parameters $\gamma _{i,j}$. Noticing that $D_{r,r}=s_{r,r}$
and using the formula 
$\varphi ^{\sharp }_q(0,r)=\prod _{l=1}^q
\frac{1}{d_{r,r+l}}$
obtained in the proof of Theorem \ref{convergenta}, we deduce that
$$\begin{array}{rcl}
\displaystyle\frac{D_{r,q}}{D_{r+1,q}}&=& 
s_{r,r}\det \Delta (\cH _{r,r},\cH _{r+1,q}) \\
 & &\\
  &=&s_{r,r}\displaystyle\frac{\prod _{r\leq k<j\leq q}d^2_{k,j}}
{\prod _{r+1\leq k\leq j\leq q}d^2_{k,j}} \\
 & & \\
 &=&s_{r,r}\prod _{j=1}^qd^2_{r,r+j}=
\displaystyle\frac{s_{r,r}}{|\varphi ^{\sharp }_q(0,r)|^2},
\end{array}
$$
which is \eqref{3?}. This relation and Theorem~ 6.2.2 in \cite{Co}
imply \eqref{3??}.

Using again \eqref{formula}, we deduce for $n<r$, that  
$$\begin{array}{rcl}
\det \Delta (\cH _{0,n},\cH _{n+1,r})&=&
\displaystyle\frac{D_{0,r}}{D_{0,n}D_{n+1,r}} \\
  & & \\
 &=&\prod _{0\leq k\leq n<j\leq r}d^2_{k,j},
\end{array}
$$
hence 
$$\det \Delta (\cH _{0,n},\cH _{n+1,\infty })
=\lim _{r\rightarrow \infty }
\det \Delta (\cH _{0,n},\cH _{n+1,r})=
\prod _{0\leq k\leq n<j}d^2_{k,j}.$$

On the other hand, 
$$\begin{array}{rcl}
\displaystyle\frac{\prod _{l=0}^ng_l}{D_{0,n}} & =& 
\displaystyle\frac{\prod _{l=0}^n\prod _{j\geq 1}d^2_{l,l+j}}
{\prod _{0\leq k<j\leq n}d^2_{k,j}} \\
 & & \\
& =&\prod _{0\leq k\leq n<j}d^2_{k,j},
\end{array}
$$ 
which shows that 
$$\det \Delta (\cH _{0,n},\cH _{n+1,\infty })=
\frac{\prod _{l=0}^ng_l}{D_{0,n}},
$$
hence \eqref{3???}.
\end{proof}

Formula \eqref{3??} would represent an analogue of \eqref{trei}, 
while \eqref{3???} is an analogue of \eqref{patru}. 
It would be of interest to express the limit
in \eqref{3???} in terms of the spectral factor $\Theta _{\phi }$.

\section{Some examples}
\subsection{Polynomials on the unit circle}
Consider $\cA =\{1-\cJ (X_1)X_1\}$. 
In this case the index set is $\NN _0$ and if $\phi $ is a linear functional 
on $\cT _1(\cA )$, then 
$$K_{\phi }(n+k,m+k)=K_{\phi }(n,m), \quad m,n,k\in \NN _0,$$
which means that $K_{\phi }$ is a Toeplitz kernel. It turns out that 
the parameters $\{\gamma _{k,j}\}$ also satisfy the Toeplitz condition, 
$\gamma _{n+k,m+k}=\gamma _{n,m}$, $n<m$, $k\geq 1$. The orthogonal 
polynomials associated to $\phi $ are then the orthogonal polynomials on the 
unit circle and \eqref{primarelatie}, \eqref{adouarelatie}  
reduce to the classical recurrence equations in \cite{Sz}. Also, 
Theorem ~\ref{convergenta} and Theorem ~\ref{szego} 
reduce to the classical results of Szeg\"o, \cite{Sz}.

\subsection{Polynomials on the real line}
Consider $\cA =\{X_1-\cJ (X_1)\}$. 
In this case the index set is still $\NN _0$, 
and if $\phi $ is a linear functional 
on $\cT _1(\cA )$, this time the kernel $K_{\phi }$
has the Hankel property, that is
$$K_{\phi }(n,m+k)=K_{\phi }(n+k,m), \quad m,n,k\in \NN _0.$$
The parameters
$\{\gamma _{k,j}\}$ do not necessarely satisfy a similar Hankel property.
In fact, it might be interesting to find a characterization
of those families of parameters $\{\gamma _{k,j}\}$ producing Hankel forms. 
Traditionally, there are other parameters, usually called 
canonical moments, that are used. The canonical moments of $\phi $ can be 
calculated by using a $Q$-$D$ (quotient-difference) algorithm (see \cite{He}).
Also, the recurrence formulas of type \eqref{primarelatie},
\eqref{adouarelatie} are replaced by a three-term recurrence equation,
\begin{equation}\label{treitermeni}
x\varphi _n(x)=b_{n}\varphi _{n+1}(x)+
a_{n}\varphi _n(x)+
b_{n-1}\varphi _{n-1}(x),
\end{equation}
with initial conditions $\varphi _{-1}=0$, $\varphi _0=1$ (\cite{Sz}).

Still, parameters $\{\gamma _{k,j}\}$ can be associated 
such that \eqref{primarelatie}, \eqref{adouarelatie} hold. Also, 
Theorem ~\ref{convergenta} and Theorem ~\ref{szego} provide
asymptotic properties of the orthogonal polynomials and, respectively, 
Hankel determinants in the corresponding Szeg\"o class.

We consider an example computing the parameters $\gamma _{k,j}$
of the Hilbert matrix,
$$H=\left[\begin{array}{cccc}
1 & \frac{1}{2} & \frac{1}{3} & \ldots  \\
  & & &  \\
\frac{1}{2} & \frac{1}{3} & \frac{1}{4} & \ldots \\
  & & &  \\
\frac{1}{3} & \frac{1}{4} & \frac{1}{5} & \ldots \\
 \vdots & \vdots & \vdots & \ddots 
\end{array}\right]
$$
and notice that the associated orthogonal polynomials
satisfy the three-term recurrence equation
\eqref{treitermeni}
with 
$$b_{n-1}=\frac{n}{2\sqrt{4n^2-1}},\quad n\geq 1,$$
and 
$$a_n=\frac{1}{2},\quad n\geq 0.$$
For example, the first $5$ polynomials are:
$$\varphi _{-1}=0, \quad \varphi _0=1, \quad \varphi _1(x)=\sqrt{3}(2x-1),$$
$$\varphi _2(x)=\sqrt{5}(6x^2-6x+1),\quad 
\varphi _3(x)=\sqrt{7}(20x^3-30x^2+12x-1).$$
The canonical moments $\{p_n\}_{n\geq 0}$ can be calculated from the 
continued fraction expansion of the Stieltjes transform of the 
uniform measure on $[0,1]$, 
$$\int _0^1\frac{dx}{z-x}=
\displaystyle\frac{1}{z-\displaystyle\frac{\frac{1}{2}}
{1-\displaystyle\frac{\frac{2}{3}\frac{1}{3}}
{z-\displaystyle\frac{\frac{2}{3}\frac{1}{2}}{1-\ldots }}}},
$$
which gives
$$p_{2k-1}=\frac{1}{2}, \quad p_{2k}=\frac{k}{2k+1}, \quad k\geq 1.
$$
We deduce that, for $n\geq 1$, 
$$\det \left[
\begin{array}{cccc}
1 & \frac{1}{2} & \ldots & \frac{1}{n}   \\
 & & & \\
\frac{1}{2} & \frac{1}{3} &  \ldots  & \frac{1}{n+1} \\
 \vdots & \vdots &  \ddots & \\
 \frac{1}{n}  & \frac{1}{n+1} & \ldots &  \frac{1}{2n} 
\end{array}
\right]
=\left(\prod _{k=1}^{n+1}\frac{1}{2k-1}\right)
\prod _{l=0}^{n-1}\prod _{k=1}^{n-l}\left(\frac{k}{k+2l+1}\right)^2.
$$
This formula, \eqref{primarelatie}, and \eqref{treitermeni} give that
$$\gamma _{0,l}=(-1)^{l-1}\frac{\sqrt{2l+1}}{l+1}, \quad l\geq 1.$$
Extending this argument (based on results from \cite{Sz}), we deduce
that
$$\gamma _{k,k+l}=(-1)^{l-1}\frac{\sqrt{(2k+1)(2k+2l+1)}}{2k+l+1},\quad 
k\in \NN _0, l\geq 1,$$
hence
$$d_{k,k+l}=\frac{l}{2k+l+1}.
$$
These formulae show that the uniform measure on $[0,1]$ does not
belong to the Szeg\"o class.

We can obtain explicit computation of $\{\gamma _{k,j}\}$
for other classes of classical orthogonal polynomials. The main point is 
to notice that if $\{\varphi \}_{n\geq 0}$ is the family of orthogonal
polynomials associated to a certain weight $w(x)$, then 
$\{\varphi _n(x,l)\}_{n\geq 0}$ is the family of orthogonal
polynomials associated to the weight $x^{2l}w(x)$. The polynomials
associated to $x^{2l}w(x)$
are called the modified orthogonal polynomials and their calculation 
for Hermite and Gegenbauer polynomials can be found, for instance, 
in \cite{DX}.
Then Theorem ~\ref{nouaparam} can be used to determine the parameters
$\{\gamma _{k,j}\}$. Details will appear in \cite{Ba}.

\subsection{Szeg\"o polynomials on several non-commuting variables}
The next examples are motivated in part by multiscale processes.
These are stochastic processes indexed by the nodes of a tree. 
Isotropic
processes on homogeneous trees were systematically studied, see
\cite{BBCGNW} and the references therein. An extension to chordal graphs was
recently given in \cite{FLW}. Some classes of stochastic 
processes associated to the full binary (Cayley) tree were
also considered (see, for instance, \cite{Fr}).

We discuss here this last example. The vertices of the Cayley tree are
indexed by $\FF ^+_N$. Let $(X,\cM ,\mu )$ be a probability space
and let $\{v_{\sigma }\}_{\sigma \in \FF ^+_N}\subset L^2(\mu )$
be a family of random variables.  Its covariance kernel is
$$K(\sigma , \tau )=\int _X\overline{v}_{\sigma }v_{\tau }dP.$$
The processes is called stationary
(see \cite{Fr}), if 
\begin{equation}\label{sta1}
K(\tau \sigma ,\tau \sigma ')=K(\sigma ,\sigma '), \quad 
\tau ,\sigma ,\sigma '\in \FF _N^+,
\end{equation}
\begin{equation}\label{sta2}
K(\sigma ,\tau )=0 \quad \mbox{ if there is no $\alpha \in \FF _N^+$
such that $\sigma =\alpha \tau $ or $\tau =\alpha \sigma $}.
\end{equation}

Let $\cA _S=\{1-\cJ (X_k)X_k\mid k=1,\ldots ,N\}\cup
\{\cJ (X_k)X_l, k,l=1,\ldots ,N, k\ne l\}$
and note that the index set of $\cA _S$ is $\FF ^+_N$. We see that 
$\phi $ is a positive functional on $\cT _N(\cA _S)$ 
if and only if $K_{\phi }$ 
is the covariance of a stationary process as above. 
It was noticed in \cite{CJ} that this happens if and only if
\begin{equation}\label{sta3}
\gamma _{\tau \sigma ,\tau \sigma '}=\gamma _{\sigma ,\sigma '}, \quad 
\tau ,\sigma ,\sigma '\in \FF _N^+,
\end{equation}
\begin{equation}\label{sta4}
\gamma _{\sigma ,\tau }=0 \quad \mbox{ if there is no $\alpha \in \FF _N^+$
such that $\sigma =\alpha \tau $ or $\tau =\alpha \sigma $},
\end{equation}
where $\{\gamma _{\sigma , \tau }\mid \sigma ,\tau \in \FF _N^+, 
\sigma \preceq \tau \}$
is the family of parameters associated to $K_{\phi }$ by Theorem~1.5.3 
in \cite{Co}.
The main consequence of these relations is that we can 
define $\gamma _{\sigma }=
\gamma _{\emptyset ,\sigma }$, $\sigma \in \FF _N^+$,
and 
$\{\gamma _{\sigma , \tau }\mid \sigma ,\tau \in \FF _N^+, \sigma \preceq
\tau \}$
is uniquely determined by $\{\gamma _{\sigma }\}_{\sigma \in \FF _N^+}$
due to the relation 
\begin{equation}\label{explo}
\left[\gamma _{\sigma ,\tau }
\right]_{|\sigma |=j,|\tau |=k}=
(\left[\gamma _{\sigma ',\tau '}
\right]_{|\sigma '|=j-1,|\tau '|=k-1})^{\oplus N}, \quad j,k\geq 1.
\end{equation}
From now on we assume that $\phi $ is unital, $\phi (1)=1$.
Then we can show that the recurrence equations \eqref{primarelatie} and 
\eqref{adouarelatie} simplify to 
$\varphi _{\emptyset }=1$ and for $k\in \{1,\ldots ,N\}$,
$\sigma \in \FF _N^+$, 
\begin{equation}\label{bszego}
\varphi _{k\sigma }=\frac{1}{d_{k\sigma }}
(X_k\varphi _{\sigma }-\gamma _{k\sigma } 
\varphi ^{\sharp}_{k\sigma -1}),
\end{equation}
where $\varphi ^{\sharp}_{\emptyset }=1$ and for $k\in \{1,\ldots ,N\}$,
$\sigma \in \FF _N^+$,
\begin{equation}\label{bsarp}
\varphi ^{\sharp}_{k\sigma }=\frac{1}{d_{k\sigma }}
(-\overline{\gamma }_{k\sigma }X_k\varphi _{\sigma }+ 
\varphi ^{\sharp}_{k\sigma -1}).
\end{equation}

We also notice that the algebra $\cT _N$ is naturally embedded
into $\cR _1$ and $\Phi _n, \Phi ^{\sharp }_n\in \cT _N$. Then, 
Theorem ~3.3 implies that $\Theta _{\phi }$ belongs to $\cT _N$
(but this is also seen directly), and through the isomorphisms 
mentioned in Section ~2.1, $\Theta _{\phi }$ can be identified 
with an element of the full Fock space over $\CC ^N$, and therefore
with a formal power series on variables $X_1$, $\ldots $, $X_N$.
Similarly, $(\Phi ^{\sharp }_n)^{-1}$ can be identified
with a formal power series, denoted $(\varphi ^{\sharp }_n)^{-1}$, 
on variables $X_1$, $\ldots $, $X_N$. 
Finally, we notice that $\phi $ belongs to the Szeg\"o class
if and only if 
$$\prod _{\sigma \in \FF ^+_N}d_{\sigma }>0.$$
For two formal series on variables $X_1$, $\ldots $, $X_N$, 
the sign $\rightarrow $ means coefficient-wise convergence.
The next result is a consequence of Theorem ~\ref{convergenta}.

\begin{theorem}\label{main}
Let $\phi $ be a functional on  
$\cT _N(\cA _S)$ and belonging to the Szeg\"o class. Then 
\begin{equation}\label{mainunu}
\phi _n\rightarrow 0
\end{equation}
and 
\begin{equation}\label{maindoi}
(\phi _n^{\sharp })^{-1}\rightarrow \Theta _{\phi }.
\end{equation}
\end{theorem}

As a consequence of Theorem ~\ref{szego} we obtain the following
result (which, aside the new geometrical interpretation, whould be also a 
direct consequence of Theorem 6.4.5 in \cite{Co}).

\begin{theorem}\label{cszego}
Let $\phi $ be a functional on  
$\cT _N(\cA _S)$ and belonging to the Szeg\"o class. Then 
\begin{equation}\label{5}
\lim _{|\tau |\rightarrow \infty }
\frac{D_{\emptyset ,\tau}}{D_{1,\tau}}=
|\Theta _{\phi }(0)|^2=\prod _{\sigma \in \FF ^+_N}d^2_{\sigma }.
\end{equation}
If we denote the above limit by $g$ and  
$$L=\prod _{\sigma \in \FF ^+_N}
d_{\sigma }^{2|\sigma |}>0,$$
then 
\begin{equation}\label{6}
\lim _{|\tau |\rightarrow \infty }
\frac{D_{\emptyset ,\tau}}{g^{l(\tau )}}
=\frac{1}{L}.
\end{equation}
\end{theorem}

Finally, we mention that similar results can be obtained for the 
commutative case, $\cA_C=\cA _S\cup \{X_kX_l-X_lX_k\mid k,l=1,\ldots ,N\}$.
Positive functionals on $\cA _C$ correspond to positive definite functions
on $\ZZ^N$. Details as well as other examples will be given in \cite{Ba}.

\end{document}